\documentclass[runningheads,a4paper]{llncs}

\usepackage{amssymb}
\usepackage{graphicx,longtable,wrapfig}
\usepackage{xcolor}
\usepackage{amsmath}
\usepackage{epstopdf}

\usepackage{url}

\usepackage{caption,subcaption}
\begin{document}
\newcommand{\E}{\mathsf{E}}
\newcommand{\Prob}{\mathsf{P}}
\newcommand{\G}{\mathsf{\Gamma}}
\newcommand{\BA}{Barab\'asi--Albert}
\newcommand{\BO}{Buckley--Osthus}
\newcommand{\HK}{Holme and Kim}

\newcommand{\widesim}[2][1.5]{
  \mathrel{\overset{#2}{\scalebox{#1}[1]{$\sim$}}}
}

\mainmatter

\title{Local clustering coefficient in generalized preferential attachment models\thanks{†This work is supported by the Russian President grant MK-527.2017.1.
} \thanks{This is an extended version of the paper appeared in Proc. WAW'15, LNCS 9479, pp. 15-28, 2015.} }

\titlerunning{Local clustering coefficient in generalized preferential attachment models}

\author{Alexander~Krot\inst{1}
\and Liudmila~Ostroumova~Prokhorenkova\inst{1,2}}

\institute{Moscow Institute of Physics and Technology, Moscow, Russia
\and
Yandex, Moscow, Russia}

\maketitle

\begin{abstract}
In this paper, we analyze the local clustering coefficient of preferential attachment models. A general approach to preferential attachment was introduced in~\cite{GPA}, where a wide class of models (PA-class) was defined in terms of constraints that are sufficient for the study of the degree distribution and the clustering coefficient. It was previously shown that the degree distribution in all models of the PA-class follows a power law. Also, the global clustering coefficient was analyzed and a lower bound for the average local clustering coefficient was obtained.
We expand the results of~\cite{GPA} by analyzing the local clustering coefficient for the PA-class of models. Namely, we analyze the behavior of $C(d)$ which is the average local clustering for the vertices of degree $d$.

\vspace{5pt}

\noindent{\bf Keywords:} networks, random graph models, preferential attachment, clustering coefficient.

\end{abstract}

\section{Introduction}\label{Introduction}

Nowadays there are a lot of practical problems connected with the analysis of growing real-world networks, from Internet and society networks~\cite{BA_Review,Networks,Chayes} to biological networks~\cite{BioInfoPrior}. Models of real-world networks are used in physics, information retrieval, data mining, bioinformatics, etc. An extensive review of real-world networks and their applications can be found elsewhere (e.g., see~\cite{BA_Review,Networks,Math_Results,Lescovec}).

It turns out that many real-world networks of diverse nature have some typical properties:
small diameter, power-law degree distribution, high clustering, and others \cite{Girvan,Newman2,Newman3,Watts}. Probably the most extensively studied property of networks is their vertex degree distribution. For the majority of studied real-world networks, the portion of vertices with degree $d$ was observed to decrease as $d^{-\gamma}$, usually with $2 < \gamma < 3$~\cite{BA1,BA2,Networks,Broder,F-F-F}.

Another important characteristic of a network is its clustering coefficient, which has the following two most used versions: the global clustering coefficient and the average local clustering coefficient (see Section~\ref{ClusteringCoefficient} for the definitions).
It is believed that for many real-world networks both the average local and the global  clustering coefficients tend to non-zero limit as the network becomes large.
Indeed, in many observed networks the values of both clustering coefficients are considerably high~\cite{Newman3}.

The most well-known approach to modeling complex networks is the prefer\-ential-attachment idea.
Many different models are based on this idea: LCD~\cite{LCD_degrees}, Buckley-Osthus~\cite{Buckley_Osthus}, Holme-Kim~\cite{Holme_Kim}, RAN~\cite{RAN}, and many others.
A general approach to preferential attachment was introduced in \cite{GPA}, where a wide class of models was defined in terms of constraints that are sufficient for the study of the degree distribution (PA-class) and the clustering coefficient (T-subclass of PA-class).

In this paper, we analyze the behavior of $C(d)$ --- the average local clustering coefficient for the vertices of degree $d$ --- in the T-subclass.
It was previously shown that in real-world networks $C(d)$ usually decreases as $d^{-\psi}$ with some parameter $\psi>0$ \cite{Catanzaro,Serrano,Vazquez}.
For some networks, $C(d)$ scales as a power law $C(d) \sim d^{-1}$ \cite{Lescovec,Ravasz}.
In the current paper, we prove that in \textit{all} models of the T-subclass the local clustering coefficient $C(d)$ asymptotically behaves as $C \cdot d^{-1}$, where $C$ is some constant. We also illustrated these results empirically. In addition, we suggested and empirically verified (for $A \le 0.75$) an approximation for the average local clustering coefficient $C_2(n)$.

The remainder of the paper is organized as follows.
In Section~\ref{PAclass}, we give a formal definition of the PA-class and
present some known results.
Then, in Section~\ref{PAlocal}, we state new results on the behavior of local clustering $C(d)$.
We prove the theorems in Section~\ref{Proves}. In Section~\ref{experiments} we make some simulations in order to illustrate our results for $C(d)$ and to empirically analyze the local clustering coefficient. Section~\ref{Conclusion} concludes the paper.

\section{Generalized Preferential Attachment}\label{PAclass}

\subsection{Definition of the PA-class}\label{class}

In this section, we define the PA-class of models which was first suggested in \cite{GPA}.
Let $G_{m}^n$ ($n \ge n_0$) be a graph with $n$ vertices $\{1, \ldots, n\}$ and $mn$ edges obtained as a result of the following process. We start at the time $n_0$ from an arbitrary graph $G_{m}^{n_0}$ with $n_0$ vertices and $m n_0$ edges. On the $(n+1)$-th step ($n\geq n_0$), we make the graph $G_{m}^{n+1}$ from $G_{m}^{n}$ by adding a new vertex $n+1$ and $m$ edges connecting this vertex to some $m$ vertices from the set $\{1, \ldots , n, n+1\}$. Denote by $d_v^{n}$ the degree of a vertex $v$ in~$G_m^n$. If for some constants $A$ and $B$ the following conditions are satisfied
\begin{equation}\label{OneStepChangedDegreeDraft}
\Prob\left( d_v^{n+1} = d_v^{n} \mid G_m^{n}\right) = 1 - A \frac{d_v^{n}}{n} - B\frac{1}{n} + O\left(\frac{\left(d_v^n\right)^2}{n^2}\right),\,\,1 \le v \le n \;,
\end{equation}
\begin{equation}\label{OneStepChangedDegreeDraft_2}
\Prob\left( d_v^{n+1} = d_v^{n} + 1 \mid G_m^{n}\right) =  A \frac{d_v^{n}}{n} + B\frac{1}{n} + O\left(\frac{\left(d_v^n\right)^2}{n^2}\right), \,\,1 \le v \le n \;,
\end{equation}
\begin{equation}\label{OneStepChangedDegreeDraft_3}
\Prob\left( d_v^{n+1} = d_v^{n} + j \mid G_m^{n}\right) =  O\left(\frac{\left(d_v^n\right)^2}{n^2}\right), \,\,
2\le j \le m,\,\,1 \le v \le n \;,
\end{equation}
\begin{equation}\label{LoopProbability}
\Prob( d_{n+1}^{n+1} =  m + j ) = O\left(\frac 1 n \right), \,\,
1\le j \le m\;,
\end{equation}
then the random graph process $G_m^n$ is a model from the PA-class.
Here, as in~\cite{GPA}, we require $2mA + B = m$ and $0 \le A \le 1$.

As it is explained in \cite{GPA}, even fixing values of parameters $A$ and $m$ does not specify a concrete procedure for constructing a network.
There are a lot of models possessing very different properties and satisfying the conditions~(\ref{OneStepChangedDegreeDraft}--\ref{LoopProbability}), e.g., the LCD, the Buckley--Osthus, the Holme--Kim, and the RAN models.

\subsection{Power Law Degree Distribution}\label{DegreeDistribution}

Let $N_n(d)$ be the number of vertices of degree $d$ in~$G_m^n$.
The following theorems on the expectation of $N_n(d)$ and its concentration were proved in~\cite{GPA}.

\begin{theorem}\label{Expectation}
For every model in PA-class and for every $d \ge m$
\vspace{-3pt}
$$
\E N_n(d) = c(m,d) \left(n + O\left(d^{2 + \frac{1}{A}}\right)\right)\,,
\vspace{-5pt}
$$
where
\vspace{-3pt}
\begin{equation*}\label{Constant}
c(m,d) = \frac{\G\left(d + \frac{B}{A}\right)\G\left(m + \frac{B+1}{A}\right)}{A\,\G\left(d + \frac{B+A+1}{A}\right)\G\left(m + \frac{B}{A}\right)}
\widesim{d \rightarrow \infty} \frac{\G\left(m + \frac{B+1}{A}\right)d^{-1-\frac{1}{A}}}{A\, \G\left(m + \frac{B}{A}\right)}
\vspace{-5pt}
\end{equation*}
and $\G(x)$ is the gamma function.
\end{theorem}

\begin{theorem}
\label{Concentration}
For every model from the PA-class and for every $d=d(n)$ we have
$$
\Prob\left(|N_n(d) - \E N_n(d)| \ge d \, \sqrt{n}  \, \log{n}\right) ={ n^{-\Omega(\log{n})}}.
$$
Therefore, for any $\delta > 0$ there exists a function $\varphi(n) \in o(1)$ such that
$$
\lim_{n \to \infty} \Prob \left(\exists \, d \le n^{\frac{A-\delta}{4A+2}}:|N_n(d) - \E N_n(d)| \ge \varphi(n)\,\E N_n(d)\right) = 0 \;.
$$
\end{theorem}
These two theorems mean that the degree distribution follows
(asymptotically) the power law with the parameter $1+\frac{1}{A}$.

\subsection{Clustering Coefficient}\label{ClusteringCoefficient}

A T-subclass of the PA-class was introduced in~\cite{GPA}.
In this case, the following additional condition is required:
\vspace{-3pt}
\begin{equation}\label{D_definition}
\Prob\left( d_i^{n+1} = d_i^{n} + 1, d_j^{n+1} = d_j^{n} + 1 \mid G_m^{n}\right) = e_{ij} \frac{D}{mn} + O\left(\frac{d_i^{n} d_j^{n}}{n^2}\right) \;.
\vspace{-3pt}
\end{equation}
Here $e_{ij}$ is the number of edges between vertices $i$ and $j$ in $G_m^n$ and $D$ is a positive constant. Note that this property still does not define the correlation between edges completely, but it is sufficient for studying both global and average local clustering coefficients.

Let us now define the clustering coefficients.
The \emph{global clustering coefficient} $C_1(G)$ is the ratio of three times the number of triangles to the number of pairs of adjacent edges in $G$. The \emph{average local clustering coefficient} is defined as follows: $C_2(G) = \frac{1}{n} \sum_{i=1}^n C(i)$, where $C(i)$ is the local clustering coefficient for a vertex $i$: $C(i) = \frac{T^i}{P_2^i}$, where $T^i$ is the number of edges between neighbors of the vertex $i$ and $P_2^i$ is the number of pairs of neighbors.
Note that both clustering coefficients are defined for graphs without multiple edges.

The following theorem on the global clustering coefficient in the T-subclass was proven in~\cite{GPA}.


\begin{theorem}\label{Cluster}
Let $G_m^n$ belong to the T-subclass with $D>0$. Then, for any $\varepsilon>0$
\begin{itemize}
\item[(1)]
If $2A<1$, then \textbf{whp}
$ \frac{6(1-2A)D - \varepsilon}{m(4(A+B) +m-1)} \le C_1(G_m^n) \le  \frac{6(1-2A)D + \varepsilon}{m(4(A+B)+m-1)} \;;$
\item[(2)]
If $2A=1$, then \textbf{whp}
$\frac{6 D - \varepsilon }{m(4(A+B)+m-1) \log n} \le C_1(G_m^n) \le \frac{6 D + \varepsilon }{m(4(A+B)+m-1) \log n} \;;$ 
\item[(3)]
If $2A>1$, then \textbf{whp}
$n^{1-2A-\varepsilon} \le C_1(G_m^n) \le n^{1-2A+\varepsilon}\;.$
\end{itemize}
\end{theorem}

Theorem \ref{Cluster} shows that in some cases ($2A\ge1$) the global clustering coefficient $C_1(G_m^n)$ tends to zero as the number of vertices grows.

The average local clustering coefficient $C_2(G_m^n)$ was not fully analyzed previously, but it was shown in \cite{GPA} that $C_2(G_m^n)$ does not tend to zero for the T-subclass with $D>0$.
In the next section, we fully analyze the behavior of the average local clustering coefficient for the vertices of degree $d$.

\section{The Average Local Clustering for the Vertices of Degree $d$}\label{PAlocal}

In this section, we analyze the asymptotic behavior of $C(d)$~--- the average local clustering for the vertices of degree $d$.
Let $T_n(d)$ be the number of triangles on the vertices of degree $d$ in $G_m^n$
(i.e., the number of edges between the neighbors of the vertices of degree $d$).
Then, $C(d)$ is defined in the following way:
\vspace{-3pt}
\begin{equation}\label{CD_definition}
C(d)= \frac{T_n(d)}{N_n(d){d \choose 2}} \,.
\vspace{-3pt}
\end{equation}
In other words, $C(d)$ is the local clustering coefficient averaged over all vertices of degree $d$.
In order to estimate $C(d)$ we should first estimate $T_n(d)$.
After that, we can use Theorems~\ref{Expectation} and~\ref{Concentration} on the behavior of $N_n(d)$.

We prove the following result on the expectation of $T_n(d)$.

\begin{theorem}\label{PATexp}
Let $G_m^n$ belong to the T-subclass of the PA-class with $D>0$. Then
\begin{itemize}
\item[(1)] if $2A<1$, then
$\E T_n(d) = K(d)\left(n + O \left( d^{2+\frac{1}{A}}\right)\right)$;
\item[(2)] if $2A = 1$, then
$\E T_n(d) = K(d)\left(n + O \left( d^{2+\frac{1}{A}} \cdot \log(n) \right)\right)$;
\item[(3)] if $2A > 1$, then
$\E T_n(d) = K(d)\left(n + O \left( d^{2+\frac{1}{A}} \cdot n^{2A-1}\right)\right)$;
\end{itemize}
where $K(d) = c(m,d) \left(D+ \frac{D}{m} \cdot  \sum_{i=m}^{d-1} \frac{i}{Ai+B} \right)   \widesim{d \rightarrow \infty} \frac{D}{A\, m} \cdot \frac{\G\left(m + \frac{B+1}{A}\right)}{A\, \G\left(m + \frac{B}{A}\right)} \cdot d^{-\frac{1}{A}}$.
\end{theorem}

Second, we show that the number of triangles on the vertices of degree $d$ is highly concentrated around its expectation.
\begin{theorem}\label{PATconc}
Let $G_m^n$ belong to the T-subclass of the PA-class with $D>0$. Then for every $d=d(n)$
\begin{itemize}
\item[(1)] if $2A<1$: $\Prob\left(|T_n(d) - \E T_n(d)| \ge d^2 \, \sqrt{n}  \, \log{n}\right) = { n^{-\Omega(\log{n})}}$;
\item[(2)] if $2A=1$: $\Prob\left(|T_n(d) - \E T_n(d)| \ge d^2 \, \sqrt{n}  \, \log^2{n}\right) = { n^{-\Omega(\log{n})}}$;
\item[(3)] if $2A>1$: $\Prob\left(|T_n(d) - \E T_n(d)| \ge d^2 \, n^{2A-\frac{1}{2}} \, \log{n}\right) = { n^{-\Omega(\log{n})}}$.
\end{itemize}
Consequently, for any $\delta > 0$ there exists a function $\varphi(n) = o(1)$ such that
\begin{itemize}
\item[(1)] if $2A\le 1$:
$\lim_{n \to \infty} \Prob \left(\exists \, d \le n^{\frac{A-\delta}{4A+2}}:|T_n(d) - \E T_n(d)| \ge \varphi(n)\,\E T_n(d)\right) = 0$;
\item[(2)] if $2A>1$: \\ $\lim_{n \to \infty} \Prob \left(\exists \, d \le n^{\frac{A(3-4A)-\delta}{4A+2}}:|T_n(d) - \E T_n(d)| \ge \varphi(n)\,\E T_n(d)\right) = 0$.
\end{itemize}
\end{theorem}

As a consequence of Theorems~\ref{Expectation}, \ref{Concentration}, \ref{PATexp}, and \ref{PATconc}, we get the following result on the average local clustering coefficient $C(d)$ for the vertices of degree $d$ in $G_m^n$.

\begin{theorem}\label{PACD}
Let $G_m^n$ belong to the T-subclass of the PA-class. Then for any $\delta > 0$ there exists a function $\varphi(n) = o(1)$ such that
\begin{itemize}
\item[(1)] if $2A \le 1$:
$\lim_{n \to \infty} \Prob \left(\exists \, d \le n^{\frac{A-\delta}{4A+2}}:\left|C(d) - \frac{K(d)}{{d\choose 2}\,c(m,d)}\right| \ge \frac{\varphi(n)}{d}\right) = 0 $;
\item[(2)] if $2A>1$:
$\lim_{n \to \infty} \Prob \left(\exists \, d \le n^{\frac{A(3-4A)-\delta}{4A+2}}:\left|C(d) - \frac{K(d)}{{d\choose 2}\,c(m,d)}\right| \ge \frac{\varphi(n)}{d}\right) = 0 $.
\end{itemize}
Note that $\frac{K(d)}{{d\choose 2}\,c(m,d)} = \frac{2D}{d\,(d-1)\,m}\left(m +  \sum_{i=m}^{d-1} \frac{i}{Ai+B}\right) \widesim{d \rightarrow \infty} \frac{2D}{mA}\cdot d^{-1} $.
\end{theorem}

It is important to note that Theorems~\ref{PATconc} and~\ref{PACD} are informative only for $A < \frac{3}{4}$, since only in this case the value $n^{\frac{A(3-4A)-\delta}{4A+2}}$ grows. 
This restriction seems technical, i.e., one may think that more accurate estimation of error terms may fill the gap between $\frac 3 4$ and 1. However, as we discuss in Section~\ref{experiments}, it seems that for $A > \frac{3}{4}$ the error terms can make a significant contribution to $C(d)$ and the obtained asymptotic may not work. This means that it is probably impossible to estimate $C(d)$ in the whole T-subclass for $A > \frac{3}{4}$ and additional constraints are needed.


In the next section, we
first prove Theorem~\ref{PATexp}.
Then, using the Azuma--Hoeffding inequality, we prove Theorem~\ref{PATconc}. Theorem~\ref{PACD} is a corollary of Theorems~\ref{Expectation}, \ref{Concentration}, \ref{PATexp}, and \ref{PATconc}.

\vspace{-3pt}

\section{Proofs}\label{Proves}

\vspace{-3pt}

In all the proofs we use the notation $\theta(\cdot)$ for error terms.
By $\theta(X)$ we denote an arbitrary function such that $|\theta(X)| < X$.

\vspace{-3pt}

\subsection{Proof of Theorem \ref{PATexp}}

\vspace{-3pt}

We need the following auxiliary theorem.

\begin{theorem}\label{SNDlemm}
Let $W_n$ be the sum of the squares of the degrees of all vertices in a model from the PA-class. Then
\begin{itemize}
\item[(1)] if $2A<1$, then $\E W_n = O(n),$
\item[(2)] if $2A=1$, then $\E W_n = O(n \cdot \log(n)),$
\item[(3)] if $2A>1$, then $\E W_n = O(n^{2A})$.
\end{itemize}
\end{theorem}

This statement is mentioned in \cite{GPA} and it can be proved by induction. Also, let $S(n,d)$ be the sum of the degrees of all the neighbors of all vertices of degree $d$. Note that $S(n,d)$ is not greater than the sum of the degrees of the neighbors of all vertices. The last is equal to $W_n$, because each vertex of degree $d$ adds $d^2$ to the sum of the degrees of the neighbors of all vertices. So, for any $d$ we have
\vspace{-3pt}
\begin{equation}\label{equat}
\E S(n,d) \le \E W_n.
\end{equation}

Now we can prove Theorem \ref{PATexp}.
Note that we do not take into account the multiplicities of edges when we calculate the number of triangles, since the clustering coefficient is defined for graphs without multiple edges.
This does not affect the final result since the number of multiple edges is small for graphs constructed according to the model~\cite{Math_Results}.

We prove the statement of Theorem 4 by induction on $d$. Also, for each $d$ we use induction on $n$.
First, consider the case $d=m$. The expected number of triangles on any vertex $t$ of degree $m$ is equal to
$\E \sum_{(i,j) \in E(G_m^t)} {\left( e_{ij} \frac{D}{mt} + O\left(\frac{d_i^t d_j^t}{t^2}\right) \right)}$
(see \eqref{D_definition}).
As $G_m^t$ has exactly $mt$ edges, we get
$\E \sum_{(i,j) \in E(G_m^t)} {\left( e_{ij} \frac{D}{mt} + O\left(\frac{d_i^t d_j^t}{t^2}\right) \right)} \\ = D + o(1)$.
The fact that
$\E \sum_{(i,j) \in E(G_m^t)} O\left(\frac{d_i d_j}{t^2}\right) = O\left(\frac{\E W_t}{t^2}\right) = o(1)$ can be shown by induction using the conditions (1-4).
We also know (see Theorem~\ref{Expectation}) that $\E N_n(m) = c(m,m)\, n + O\left(1 \right)$. So, $\E T_n(m) = \left(D + o(1) \right) \,\left( c(m,m)\, n + O\left(1 \right)  \right) \\ =  K(m) \left(n + O\left(1 \right) \right)$. This concludes the proof for the case $d=m$ for all values of $A$ ($2A < 1$, $2A = 1$ and $2A > 1$).

Consider the case $d>m$. Note that the number of triangles on a vertex of degree $d$ is $O\left(d \right)$, since this number is $O(1)$ when this vertex appears plus at each step we get a triangle only if we hit both the vertex under consideration and a neighbor of this vertex, and our vertex degree equals $d$, therefore we get at most $d \, m$ triangles.
Also, $\E N_n(d) = c(m,d) \, \left(n + O\left(d^{2+\frac{1}{A}} \right) \right)$. So we have
$\E T_n(d) = O(d) \, c(m,d) \, \left(n + O\left(d^{2+\frac{1}{A}} \right) \right)$.
In particular, for $n \le Q \cdot d^2$ (where the constant $Q$ depends only on $A$ and $m$ and will be defined later) we have $\E T_n(d) = O\left(c(m,d) \, d^{3+\frac{1}{A}}\right) = O\left(d^2 \right) = K(d) \cdot O\left(d^{2+\frac{1}{A}}\right)$. This concludes the proof for the case $d > m$, $n \le Q d^2$ for all values of $A$.

Now, consider the case $d>m$, $n > Q \, d^2$.
Once we add a vertex $n+1$ and $m$ edges, we have the following possibilities.

1. At least one edge hits a vertex of degree $d$. Then $T_n(d)$ is decreased by the number of triangles on this vertex (because this vertex is a vertex of degree $d+1$ now). The probability to hit a vertex of degree $d$ is $\frac{A \, d + B}{n} + O\left(\frac{d^2}{n^2} \right)$. Summing over all vertices of degree $d$ we obtain that $\E T_n(d)$ is decreased by:
\begin{equation}\label{eq:1}
\left(\frac{Ad+B}{n} + O\left(\frac{d^2}{n^2} \right) \right) \cdot \E T_n(d)\,.
\end{equation}

2. Exactly one edge hits a vertex of degree $d-1$. Then $T_n(d)$ is increased by the number of triangles on this vertex. The probability to hit a vertex of degree $d-1$ once is equal to $\frac{A \, (d-1) + B}{n} + O\left(\frac{d^2}{n^2} \right)$. Summing over all vertices of degree $d-1$ we obtain that the value $\E T_n(d)$ is increased by:

\vspace{-5pt}
\begin{equation}\label{theorem2}
\left(\frac{A(d-1)+B}{n} + O\left(\frac{d^2}{n^2} \right) \right) \cdot \E T_n(d-1)\,.
\end{equation}
\vspace{-5pt}

3. Exactly one edge hits a vertex of degree $d-1$ and another edge hits its neighbor. Then, in addition to \eqref{theorem2}, $T_n(d)$ is increased by $1$. The probability to hit a vertex of degree $d-1$ and its neighbor is equal to $\frac{D}{mn} + O\left(\frac{(d-1) \, d_i}{n^2} \right)$, where $d_i$ is the degree of this neighbor. Summing over the neighbors of a given vertex of degree $d-1$ and summing then over all vertices of degree $d-1$ we obtain that $\E T_n(d)$ is increased by:
\vspace{-5pt}
\begin{multline}
(d-1)\, \E N_n(d-1)\, \frac{D}{mn} + O\left(\frac{d \cdot \E \sum_{\substack{i: i\text{ is a neighbor } \\ \text{ of a vertex of degree  $d-1$}}} d_i}{n^2} \right) \\
= (d-1)\, \E N_n(d-1)\, \frac{D}{mn} + O\left(\frac{d \, \E S(n,d)}{n^2} \right)\,.
\end{multline}
\vspace{-5pt}

4. Exactly $i$ edges hit a vertex of degree $d-i$, where $i$ is between $2$ and $m$.
If no edges hit the neighbors of this vertex, then $T_n(d)$ is increased only by the number of triangles on this vertex.
The probability to hit a vertex of degree $d-i$ exactly $i$ times is equal to $O\left(\frac{d^2}{n^2} \right)$.
If we also hit its neighbors, then $T_n(d)$ is additionally increased by $1$ for each neighbor.
The probability to hit a vertex of degree $d-i$ exactly $i$ times and hit some its neighbor is, obviously, $O\left(\frac{d^2}{n^2} \right)$. Summing over all vertices of degree $d-i$ and then summing over all $i$ from $2$ to $m$, we obtain that $\E T_n(d)$ is increased by:
\vspace{-5pt}
\begin{multline}\label{eq:4}
\sum_{i=2}^{m}{ \left(\E T_n(d-i) \cdot O\left(\frac{d^2}{n^2} \right) + O\left(\frac{d^2}{n^2} \right) \cdot (d-i) \cdot \E N_n(d-i) \right)} \\
= O\left(\frac{d^2}{n^2} \right) \, \E T_n(d) + O\left(\frac{d^3}{n^2} \right) \, \E N_n(d)\,.
\end{multline}
\vspace{-5pt}

Finally, using \eqref{eq:1}-\eqref{eq:4} and the linearity of the expectation, we get

\begin{multline}\label{Trec}
\E T_{n+1}(d) = \E T_{n}(d) - \left(\frac{Ad+B}{n} + O\left(\frac{d^2}{n^2} \right) \right) \, \E T_n(d) \\
+ \left(\frac{A(d-1)+B}{n} + O\left(\frac{d^2}{n^2} \right) \right) \, \E T_n(d-1)  + (d-1)\, \E N_n(d-1)\, \frac{D}{mn}
\\ + O\left(\frac{d \, \E S(n,d)}{n^2} \right)
+ O\left(\frac{d^2}{n^2} \right) \, \E T_n(d) + O\left(\frac{d^3}{n^2} \right) \, \E N_n(d) \\
=  \left(1-\frac{Ad+B}{n} \right)  \, \E T_n(d) + \frac{A(d-1)+B}{n} \, \E T_n(d-1)
\\ + O\left(\frac{d^2}{n^2} \right) \, \left(\E T_n(d) + \E T_n(d-1) \right) + O\left(\frac{d^3}{n^2} \right) \, \E N_n(d)
\\  + \frac{D}{mn} \, (d-1) \, \E N_n(d-1) + O\left(\frac{d\cdot \E S(n,d)}{n^2} \right)\,.
\end{multline}

Consider the case $2A < 1$ (the cases $2A = 1$ and $2A > 1$ will be analyzed similarly).
We prove by induction on $d$ and $n$ that
\begin{equation}\label{eq:induction}
\E T_n(d) = K(d) \, \left(n+\theta \left(C\cdot d^{2+\frac{1}{A}}\right) \right)
\end{equation}
for some constant $C>0$.
Let us assume that $\E T_i(\tilde d) = K(\tilde d) \, \left(i+\theta \left(C\cdot \tilde d^{2+\frac{1}{A}}\right) \right)$ for $\tilde d < d$ and all $i$ and for $\tilde d = d$
and $i < n+1$.

Recall that $K(d) = c(m,d) \left(D+ \frac{D}{m} \cdot  \sum_{i=m}^{d-1} \frac{i}{Ai+B} \right)$ and $\E N_n(d) = c(m,d) \cdot \left(n+O\left(d^{2+\frac{1}{A}}\right) \right)$. If $2A < 1$, then from (\ref{equat}) and Theorem 7 we get $\E S(n,d) = O(n)$ and we obtain:
\begin{multline*}
\E T_{n+1}(d) = \left(1-\frac{Ad+B}{n} \right) \, K(d) \, \left(n+\theta \left(C d^{2+\frac{1}{A}}\right) \right)
\\
+  \frac{A(d-1)+B}{n} \, K(d-1) \,  \left(n+\theta \left(C (d-1)^{2+\frac{1}{A}}\right) \right)
\\ + O\left(\frac{d^2}{n^2} \right) \, \left(K(d) \, \left(n+\theta \left(C d^{2+\frac{1}{A}}\right) \right)  + K(d-1) \, \left(n+\theta \left(C (d-1)^{2+\frac{1}{A}}\right) \right) \right) \\
 + O\left(\frac{d^3}{n^2} \right) \, c(m,d) \, \left(n+O\left(d^{2+\frac{1}{A}}\right) \right)
\\ + \frac{D}{mn} \, (d-1) \, c(m,d-1) \, \left(n+O\left(d^{2+\frac{1}{A}}\right) \right) + O\left(\frac{d}{n} \right)\,.
\end{multline*}

Note that $K(d) = \frac{A(d-1)+B}{Ad+B+1} \, K(d-1) + \frac{D(d-1)}{m\left(Ad+B+1 \right)} \, c(m,d-1)$. Therefore, we obtain:

\begin{multline*}
\E T_{n+1}(d) = K(d) \, (n+1) + K(d) \, \left(1-\frac{Ad+B}{n} \right) \, \theta \left(C\, d^{2+\frac{1}{A}} \right)
\\ +K(d-1) \, \frac{A(d-1)+B}{n} \,  \theta \left(C \, (d-1)^{2+\frac{1}{A}} \right)
\\ + \frac{D(d-1)}{mn} \, c(m,d) \, O\left(d^{2+\frac{1}{A}} \right) + O\left(\frac{d}{n} \right) + O\left(\frac{d^2}{n^2} \right) \, \left(K(d) \, n \right.
\\ \left. + K(d) \, \theta \left(C\, d^{2+\frac{1}{A}} \right) + K(d-1) \, n + K(d-1) \, \theta \left(C\, (d-1)^{2+\frac{1}{A}} \right) \right)
\\ + O\left(\frac{d^3}{n^2} \right) \, \left(c(m,d) \, n + c(m,d) \, O\left(d^{2+\frac{1}{A}} \right) \right)\,.
\end{multline*}

In order to show~\eqref{eq:induction}, it remains to prove that for some large enough $C$:
\begin{multline}\label{eq:remains}
K(d) \, \left(\frac{Ad+B}{n} \right) \, C \, d^{2+\frac{1}{A}} \ge  K(d-1) \,\frac{A(d-1)+B}{n} \, C \, (d-1)^{2+\frac{1}{A}} \\ + O\left(\frac{d^2}{n} \right) + O\left(C \, \frac{d^4}{n^2} \right) + O\left(\frac{d^4}{n^2} \right)\,.
\end{multline}

First, we analyze the following difference:
\begin{multline*}
K(d) \, \left(\frac{Ad+B}{n} \right) \, d^{2+\frac{1}{A}} - K(d-1) \, \frac{A(d-1)+B}{n} \,  (d-1)^{2+\frac{1}{A}} \\
 = \frac{Ad+B}{n} \, d^{2+\frac{1}{A}} \, \left(\frac{A(d-1)+B}{Ad+B+1} \, K(d-1) + \frac{D(d-1)}{m(Ad+B+1)} \, c(m,d-1) \right) \\ - \frac{A(d-1)+B}{n} \, K(d-1)  \, (d-1)^{2+\frac{1}{A}} =  \frac{(Ad+B)D(d-1)}{mn(Ad+B+1)} \, c(m,d-1) \, d^{2+\frac{1}{A}}
\\ + K(d-1)  \, \frac{A(d-1)+B}{n} \, \left(\frac{Ad+B}{Ad+B+1} d^{2+\frac{1}{A}} - (d-1)^{2+\frac{1}{A}} \right)
\\ \ge \frac{(Ad+B)D(d-1)}{mn(Ad+B+1)} \, c(m,d-1) \, d^{2+\frac{1}{A}} \\
+ (d-1)^{2+\frac{1}{A}} \, K(d-1) \, \frac{A(d-1)+B}{n} \cdot \frac{2A^2d + 2AB +B}{Ad(Ad+B+1)} \\
\ge \frac{(Ad+B)D(d-1)}{mn(Ad+B+1)} \, c(m,d-1) \, d^{2+\frac{1}{A}}\, .
\end{multline*}

Therefore, Equation \eqref{eq:remains} becomes:
$$ C \, \frac{(Ad+B)D(d-1)}{mn(Ad+B+1)} \, c(m,d-1) \, d^{2+\frac{1}{A}}
 \ge O\left(\frac{d^2}{n} \right) + O\left(C \, \frac{d^4}{n^2} \right) + O\left(\frac{d^4}{n^2} \right)\,.
$$

In the case $2A = 1$ this inequality will be:
\begin{multline*}
C \, \frac{(Ad+B)D(d-1)}{mn(Ad+B+1)} \, c(m,d-1) \, d^{2+\frac 1 A} \, \log(n) \\ \ge O\left(\frac{d^2}{n} \right) + O\left(C \, \frac{d^4 \cdot \log(n)}{n^2} \right) + O\left(\frac{d^4}{n^2} \right) + O\left(\frac{d \, \log(n)}{n} \right).
\end{multline*}

In the case $2A > 1$ this inequality will be:
\begin{multline*}
C \, \frac{(Ad+B)D(d-1)}{mn(Ad+B+1)} \, c(m,d-1) \, d^{2+\frac 1 A} \, n^{2A-1} \\ \ge O\left(\frac{d^2}{n} \right) + O\left(C \, \frac{d^4 \, n^{2A-1}}{n^2} \right) + O\left(\frac{d^4}{n^2} \right) + O\left(\frac{d \, n^{2A}}{n^2} \right).
\end{multline*}

It is easy to see that for $n \ge Q \cdot d^2$ (for some large $Q$ which depends only on the parameters of the model) these three inequalities are satisfied. This concludes the proof of the theorem.

\subsection{Proof of Theorem \ref{PATconc}}

This theorem is proved similarly to the concentration theorem from \cite{GPA}.
We also need the following notation (introduced in \cite{GPA}):

\begin{gather*}
p_n(d) = \Prob\left( d_v^{n+1} = d \mid d_v^{n} =d \right) = 1 - A \frac{d}{n} - B\frac{1}{n} + O\left(\frac{d^2}{n^2}\right) \, ,
\\
p_n^1(d) := \Prob\left( d_v^{n+1} = d+1 \mid d_v^{n} = d \right) =  A \frac{d}{n} + B\frac{1}{n} + O\left(\frac{d^2}{n^2}\right) \,,
\\
p_n^j(d) := \Prob\left( d_v^{n+1} = d+j \mid d_v^{n} = d \right) =  O\left(\frac{d^2}{n^2}\right),  \,\,
2\le j \le m \,,
\\
p_n := \sum_{k=1}^m \Prob( d_{n+1}^{n+1} =  m + k ) = O\left(\frac 1 n \right) \;.
\end{gather*}

To prove Theorem \ref{PATconc} we also need the Azuma--Hoeffding inequality:

\begin{theorem}[Azuma, Hoeffding]\label{Azuma}
Let $(X_i)_{i=0}^{n}$ be a martingale such that
$|X_i - X_{i-1}| \le c_i$ for any $1 \le i \le n$.
Then
$
\Prob\left(|X_n - X_0| \ge x \right) \le 2e^{-\frac{x^2}{2\sum_{i=1}^n c_i^2}}
$
for any $x>0$.
\end{theorem}

Consider the random variables
$X_i(d) = \E(T_n(d)\mid G_m^i)$,
$i = 0, \ldots, n$. Note that $X_0(d) = \E T_n(d)$ and $X_n(d) = T_n(d)$.
It is easy to see that $X_n(d)$ is a martingale.

We will prove below that for any $i = 0, \ldots, n-1$
\begin{itemize}
\item[(1)] if $2A<1$, then $|X_{i+1}(d) - X_{i}(d)| \le M d^2,$
\item[(2)] if $2A=1$, then $|X_{i+1}(d) - X_{i}(d)| \le M d^2 \log(n),$
\item[(2)] if $1<2A<\frac{3}{2}$, then $|X_{i+1}(d) - X_{i}(d)| \le M d^2 n^{2A-1}$,
\end{itemize}
where $M > 0$ is some constant.
The theorem follows from this statement immediately.
Indeed, consider the case $2A < 1$. Put $c_i = M d^2$ for all $i$.
Then from Azuma--Hoeffding inequality it follows that
\vspace{-0.2cm}
$$
\Prob\left(|T_n(d) - \E T_n(d)| \ge d^2 \, \sqrt{n}  \, \log{n}\right) \le
2\exp\left\{-\frac{n \,d^4 \,\log^2{n} }{2\,n\,M^2 d^4 } \right\} = { n^{-\Omega(\log{n})}}\,.
\vspace{-0.2cm}
$$
Therefore, for the case $2A < 1$ the first statement of the theorem is satisfied. If $d \le n^{\frac{A-\delta}{4A+2}}$, then the value
$n \,d^{-1/A}$ is considerably greater than $d^2 \, \log{n}\, \sqrt{n}$. From this the second statement of the theorem follows. The cases $2A = 1$ and $2A > 1$ can be considered similarly.
It remains to estimate $|X_{i+1}(d) - X_{i}(d)|$.

Fix $0 \le i \le n-1$ and some graph $G_m^{i}$.
Note that
\vspace{-0.1cm}
\begin{multline*}
\left|\E\left(T_n(d)\mid G_m^{i+1}\right) - \E\left(T_n(d)\mid G_m^{i}\right)\right| \le \\
\le \max_{\tilde G_m^{i+1}\supset G_m^{i}}  \left\{\E\left(T_n(d)\mid \tilde G_m^{i+1}\right)\right\} -
\min_{\tilde G_m^{i+1}\supset G_m^{i}}  \left\{\E\left(T_n(d)\mid \tilde G_m^{i+1}\right)\right\}.
\vspace{-0.1cm}
\end{multline*}

Put
$\hat G_m^{i+1} = \arg \max \E(T_n(d)\mid \tilde G_m^{i+1})$,
$\bar G_m^{i+1} = \arg \min \E(T_n(d)\mid \tilde G_m^{i+1})$.
It is sufficient to estimate the difference
$\E(T_n(d)\mid \hat G_m^{i+1}) - \E(T_n(d)\mid \bar G_m^{i+1})$.

For $i+1 \le t \le n$ put
\vspace{-0.2cm}
$$
\delta_t^i(d) = \E(T_t(d)\mid \hat G_m^{i+1}) - \E(T_t(d)\mid \bar G_m^{i+1}).
\vspace{-0.2cm}
$$

First, let us note that for $n \le W\cdot d^2$ (the value of constant $W$ will be defined later)
we have $\delta_n^i(d) \le \frac{2mn}{d} \cdot \left(\frac{m (m-1)}{2} + d\, m \right) \le 4m^2 n \le Md^2 \le Md^2 \log(n) \le Md^2  n^{2A-1}$ (since we have at most $\frac{2mn}{d}$ vertices of degree $d$, and each vertex of degree $d$ has at most $\frac{m(m-1)}{2}$ triangles when this vertex appears plus at each step we get a triangle only if we hit both the vertex under consideration and a neighbor of this vertex, and our vertex degree is equal to $d$, therefore we get at most $d\, m$ triangles) for some constant $M$ which depends only on $W$ and $m$.

It remains to estimate $\delta_n^i(d)$ for $n > W d^2$.
Consider the case $2A < 1$. We want to prove that $\delta^i_n(d) \le M d^2$ for $n > W d^2$ by induction.
Suppose that $n=i+1$.
Fix $G_m^{i}$.
Graphs $\hat G_m^{i+1}$ and $\bar G_m^{i+1}$ are obtained from the graph $G_m^{i}$ by
adding the vertex $i+1$ and $m$ edges.
These $m$ edges can affect the number of triangles on at most $m$ previous vertices.
For example, they can be drown to at most $m$ vertices of degree $d$ and decrease $T_i(d)$ by at most $\frac{m\,d\,(d-1)}{2}$.
Such reasonings finally lead to the estimate
$\delta^i_{i+1}(d) \le M d^2$ for some $M$.

Now let us use the induction. Consider $t$: $i+1 \le t \le n-1$, $t > W\, d^2$ (note that the smaller values of $t$ were already considered).
Using similar reasonings as in the proof of Theorem~\ref{PATexp} we get:
\vspace{-0.2cm}
\begin{equation*}
\delta^i_{t+1}(m) = \delta^i_t(m) \left(1 - p_t(m) \right) + O\left(\frac{1}{t}\right),
\vspace{-0.3cm}
\end{equation*}
\begin{multline*}\label{eq:delta}
\delta^i_{t+1}(d) = \delta^i_t(d) \left(1 - p_t(d) \right) +
\delta^i_t(d-1) \, p_t^1(d-1)
\\ + (d-1)\cdot \left(\E(N_t(d-1) \mid \hat G_m^i) - \E(N_t(d-1) \mid \bar G_m^i) \right) \cdot \frac{D}{mt}
\\+ O\left(\frac{d \cdot \E S(t,d-1)}{t^2} \right)  + O\left(\frac{\E T_t(d) \cdot d^2}{t^2} \right) + O\left(\frac{\E N_t(d)\cdot d^3}{t^2} \right)\, .
\vspace{-0.2cm}
\end{multline*}

Note that $\E (N_t(d) \mid \hat G_m^{i+1}) - \E (N_t(d) \mid \bar G_m^{i+1}) = O\left(d \right)$ (see \cite{GPA}) and $\E S(t,d-1) = O\left(t \right)$.
From this recurrent relations it is easy to obtain by induction that
$\delta^i_n(d) \le M d^2$ for some $M$. Indeed,
\vspace{-0.2cm}
\begin{multline*}
\delta^i_{t+1}(m) \le Mm^2 \left(1 - p_t(m) \right) + \frac{C_1}{t} \le Mm^2  \left(1 - \frac{Am+B}{t} + \frac{C_2}{t^2} \right) + \frac{C_1}{t} \le Mm^2
\end{multline*}
for sufficiently large $M$. By $C_i$, $i = 1, 2, \ldots$, we denote some positive constants. For $d > m$ we get
\begin{multline*}
\delta^i_{t+1}(d) \le Md^2 (1 - p_t(d)) + M(d-1)^2 p_t^1(d-1)  + C_{3}\frac{d^2}{t} + C_{4}\frac{d^4}{t^2}
\\ \le Md^2 \left(1 - \frac{Ad+B}{t} + C_5 \frac{d^2}{t^2} \right) + M(d-1)^2  \left(\frac{A(d-1)+B}{t} + C_6 \frac{d^2}{t^2} \right) + C_{3}\frac{d^2}{t}
\\ + C_{4}\frac{d^4}{t^2}
\le Md^2  + \frac{M}{t}\left(A(-3d^2 + 3d - 1) + B(-2d + 1)
+ C_7 \frac{d^4}{t} + C_3 \frac{d^2}{M} \right. \\ \left. + C_4 \frac{d^4}{Mt} \right)
 \le Md^2 + \frac{M}{t} \left(\left(-3A + C_7\frac{d^2}{t} + \frac{C_3}{M} + C_4\frac{d^2}{Mt} \right) \cdot d^2 \right. \\ + \left. \left(3A-2B\right) \cdot d + (B-A)\right) \le Md^2 \; .
\end{multline*}
for sufficiently large $W$ and $M$.

In the case $2A = 1$ we have $\E S(t,d-1) = O\left(t\log(t) \right)$ and we get the following inequalities:
\begin{equation*}
\delta^i_{t+1}(m) \le Mm^2 \log(t) \left(1 - p_t(m) \right) + \frac{C_1\log(t)}{t} \le Mm^2\log(t+1),
\end{equation*}
\begin{multline*}
\delta^i_{t+1}(d) \le Md^2 \log(t) (1 - p_t(d)) + M(d-1)^2 \log(t) \, p_t^1(d-1)
\\ + C_{2}\frac{d^2}{t} + C_{3}\frac{d \log(t)}{t} + C_{4}\frac{d^4 \log(t)}{t^2} \le Md^2 \log(t+1)\; .
\end{multline*}

In the case $2A > 1$ we have $\E S(t,d-1) = O\left(t^{2A} \right)$ and we get the following inequalities:
\begin{equation*}
\delta^i_{t+1}(m) \le Mm^2 t^{2A-1} \left(1 - p_t(m) \right) + \frac{C_1 t^{2A-1}}{t} \le Mm^2 (t+1)^{2A-1},
\end{equation*}
\begin{multline*}
\delta^i_{t+1}(d) \le Md^2 t^{2A-1}  (1 - p_t(d)) + M(d-1)^2 \, t^{2A-1} p_t^1(d-1) \\ + C_{2}\frac{d^2}{t} + C_{3}\frac{d \cdot t^{2A-1}}{t} + C_{4}\frac{d^4 t^{2A-1}}{t^2} \le Md^2 (t+1)^{2A-1}\; .
\end{multline*}

This concludes the proof of Theorem \ref{PATconc}.

\section{Experiments}\label{experiments}

In this section, we choose a three-parameter model from the family of polynomial graph models defined in \cite{GPA} and analyze the local clustering coefficients $C(d)$ and $C_2(n)$. First, we illustrate our results on $C(d)$ which we proved in the previous section. In addition, we consider the case $A>\frac{3}{4}$, for which we do not have a theoretical proof. In this case, our approximation of $C(d)$ slightly deviates from the experiment. Finally, we discuss how the average local clustering coefficient $C_2(n)$ can be approximated.

\subsection{Local Clustering Coefficient $C(d)$}

\begin{figure}[t]
\centering

\begin{subfigure}[t]{0.5\textwidth}
\centering
\includegraphics[width=\textwidth]{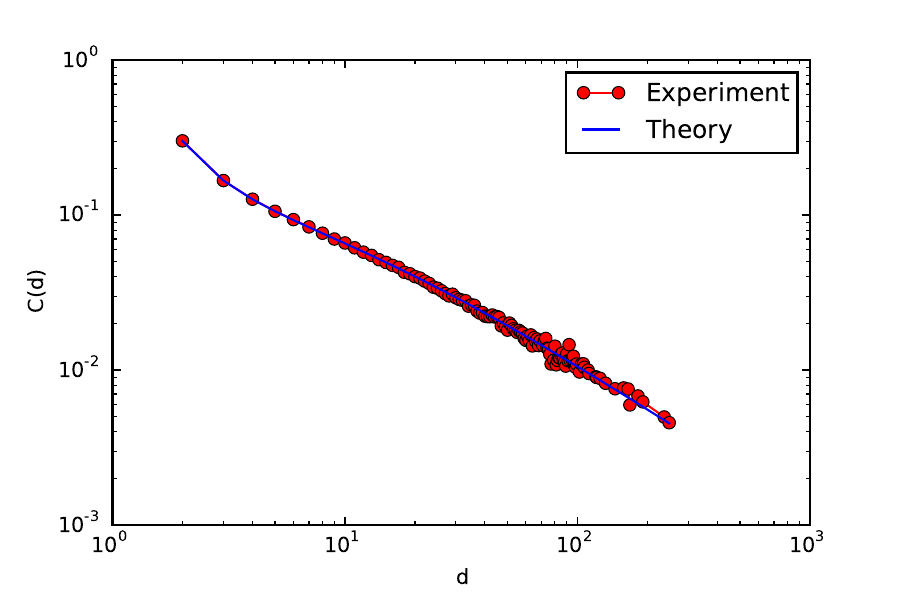}
\caption{$A=0.25$}
\label{fig:mean and std of net14}
\end{subfigure}%
\hfill
\begin{subfigure}[t]{0.5\textwidth}
\centering
\includegraphics[width=\textwidth]{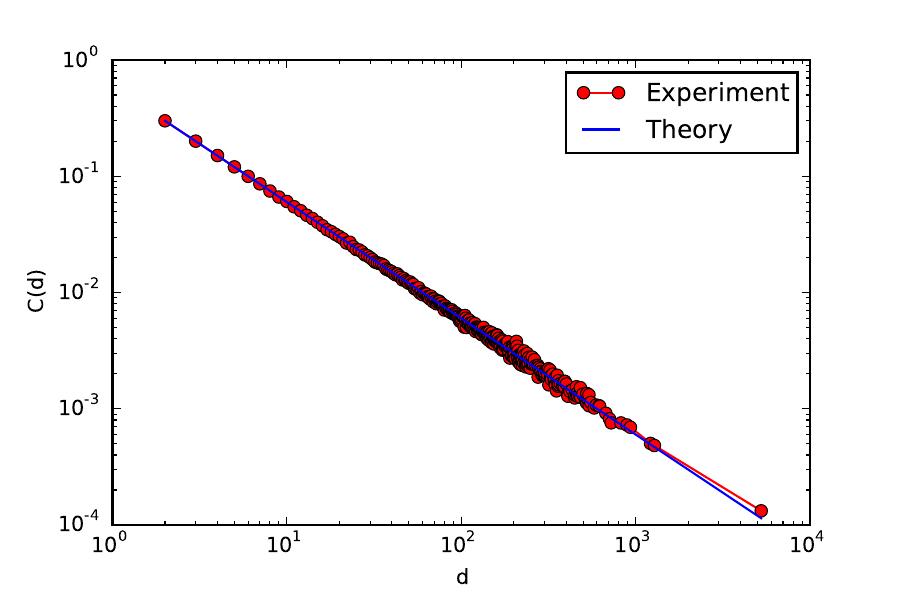}
\caption{$A=0.5$}
\label{fig:mean and std of net24}
\end{subfigure}

\bigskip

\begin{subfigure}[t]{0.5\textwidth}
\centering
\includegraphics[width=\textwidth]{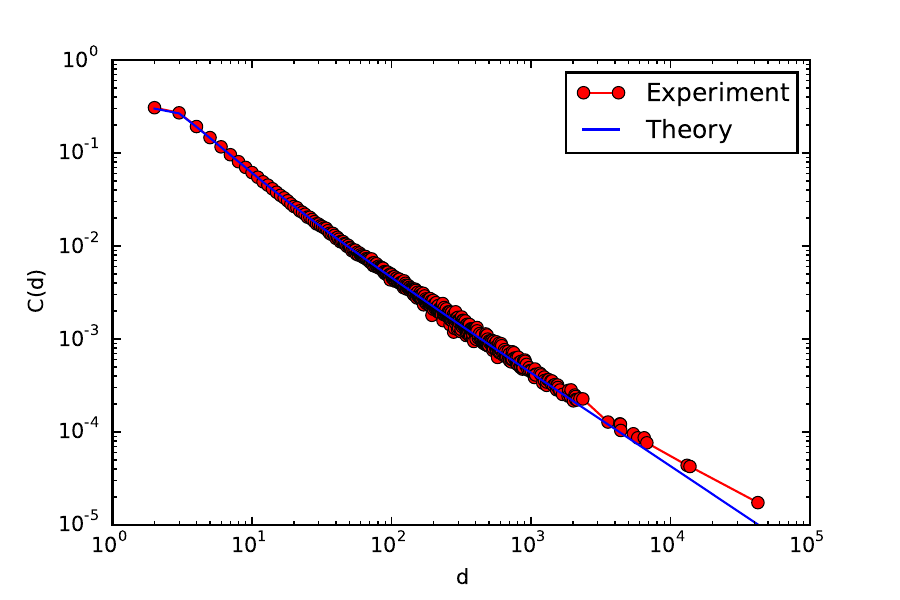}
\caption{$A=0.7$}
\label{fig:mean and std of net34}
\end{subfigure}%
\hfill
\begin{subfigure}[t]{0.5\textwidth}
\centering
\includegraphics[width=\textwidth]{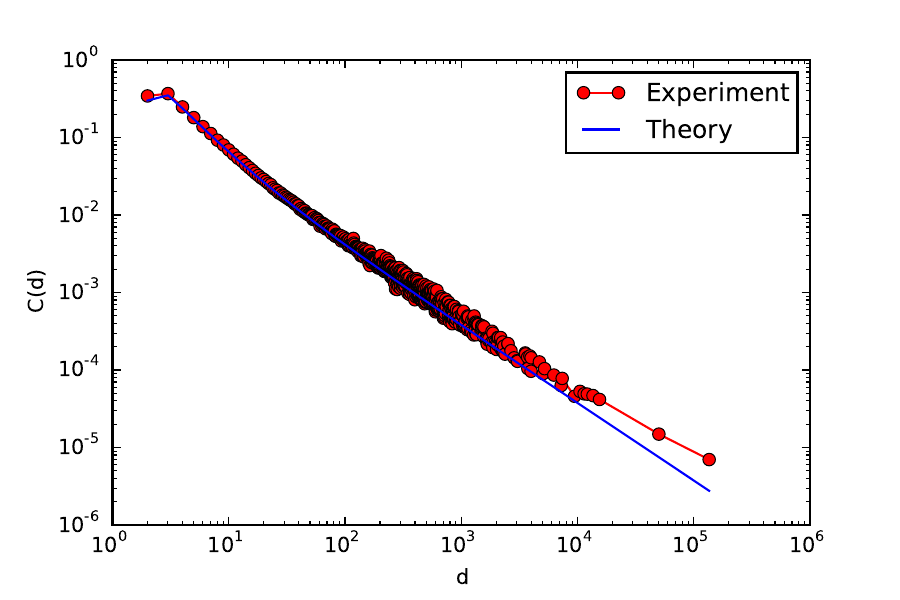}
\caption{$A=0.8$}
\label{fig:mean and std of net44}
\end{subfigure}
\caption{The behavior of $C(d)$}
\label{fig:C_D_graphs}
\end{figure}

First, we generated three polynomial graphs with $n = 10^6$, $m = 2$, $D=0.3$ and different values of $A$. In other words, we fixed the probability of a triangle formation and vary the parameter of the power-law degree distribution. Detailed graph generation process is described in \cite{GPA}. We choose $A$ to be 0.25, 0.5 and 0.7, which corresponds to the three cases of Theorems~\ref{PATexp} and~\ref{PATconc}. Also these cases correspond to three different types of a power-law degree distribution: with a finite variance, with infinite variance and the border case with $\gamma=3$. Figure~\ref{fig:C_D_graphs} illustrates our main result (Theorem~\ref{PACD}). Here the theoretical value of $C(d)$ is equal to $\frac{K(d)}{{d\choose 2}\,c(m,d)} = \frac{2D}{d\,(d-1)\,m}\left(m +  \sum_{i=m}^{d-1} \frac{i}{Ai+B}\right)$ according to Theorem~\ref{PACD}. 
We have also considered the case $A=0.8$ for which we do not have a  theoretical proof. In this case, the experimental result is also close to the theoretical approximation.
However, one can observe that our approximation slightly underestimates $C(d)$ even for small values of $d$ (see, e.g., $d=2$ on Figure~\ref{fig:mean and std of net44}).
This means that for $A > \frac 3 4$ the error terms can make a significant contribution to the value of $C(d)$ and it is probably impossible to get the accurate approximation for the whole T-subclass for such $A$. So, our restriction $A < \frac 3 4$ is essential.
In all four cases, the difference for large $d$ can be explained by the error term.

\subsection{Average Local Clustering Coefficient}

In this section we empirically analyze the average local clustering coefficient for the PA-class of models. Recall that the average local clustering coefficient is defined as: $C_2(n) = \frac{1}{n} \sum_{i=1}^{n}{C(i)}$, where $C(i)$ is the local clustering coefficient for a vertex $i$: $C(i) = \frac{T^i}{P_2^{i}}$, $T^i$ is the number of triangles on the vertex $i$ and $P_2^{i}$ is the number of pairs of neighbors. Also $C_2(n)$ can be represented in the following form: $C_2(n) = \frac{1}{n} \cdot \sum_{d=m}^{\infty} {\frac{T_n(d)}{\frac{d(d-1)}{2}}}$.

Using Theorem \ref{PATexp} we can approximate the expectation of $C_2(n)$:

\begin{multline}\label{C_2D}
\E C_2(n) = \frac{1}{n} \cdot \sum_{d=m}^{\infty} {\frac{\E T_n(d)}{\frac{d(d-1)}{2}}} = \\ = \sum_{d=m}^{\infty} { \frac{2D}{d(d-1)} \left[1 + \frac{1}{m} \sum_{i=m}^{d-1}{\frac{i}{Ai+B}} \right] \cdot} \\ \cdot{ \frac{\G\left(m+\frac{B+1}{A}\right) \G\left(d+\frac{B}{A}\right)}{A \G\left(m+\frac{B}{A}\right) \G\left(d+\frac{(B+A+1)}{A}\right)} \cdot \left[1 + O\left(\frac{X}{n} \right) \right] } = \sum_{d=m}^{\infty} {f(d) \cdot \left[1 + O\left(\frac{X}{n} \right) \right]},
\end{multline}
where:
\begin{multline*}
f(d) = \frac{2D}{d(d-1)} \left[1 + \frac{1}{m} \sum_{i=m}^{d-1}{\frac{i}{Ai+B}} \right] \cdot \frac{\G\left(m+\frac{B+1}{A}\right) \G\left(d+\frac{B}{A}\right)}{A \G\left(m+\frac{B}{A}\right) \G\left(d+\frac{(B+A+1)}{A}\right)}
\end{multline*}
and
\begin{itemize}
\item[(1)] if $2A<1$, then $X = d^{2+\frac{1}{A}}$,
\item[(2)] if $2A=1$, then $X = d^{2+\frac{1}{A}} \cdot \log(n)$,
\item[(2)] if $2A>1$, then $X = d^{2+\frac{1}{A}} \cdot \log(n^{2A})$.
\end{itemize}

It is hard to compute $\sum_{d=m}^{\infty} {f(d)}$ analytically. Moreover, it is impossible to prove that the error term $\sum_{d=m}^{\infty} {f(d) \cdot O\left(\frac{X}{n} \right)}$ behaves as $o(1)$, since this series does not converge. Therefore, in this section we empirically analyze how well $\sum_{d=m}^{\infty} {f(d)}$ approximates the clustering coefficient $C_2(n)$. Further in this section we consider the behavior of $C_2(n)$ depending on $A$ and on $D$. 

\begin{figure}[t]
\begin{center}
\includegraphics[height = 7cm]{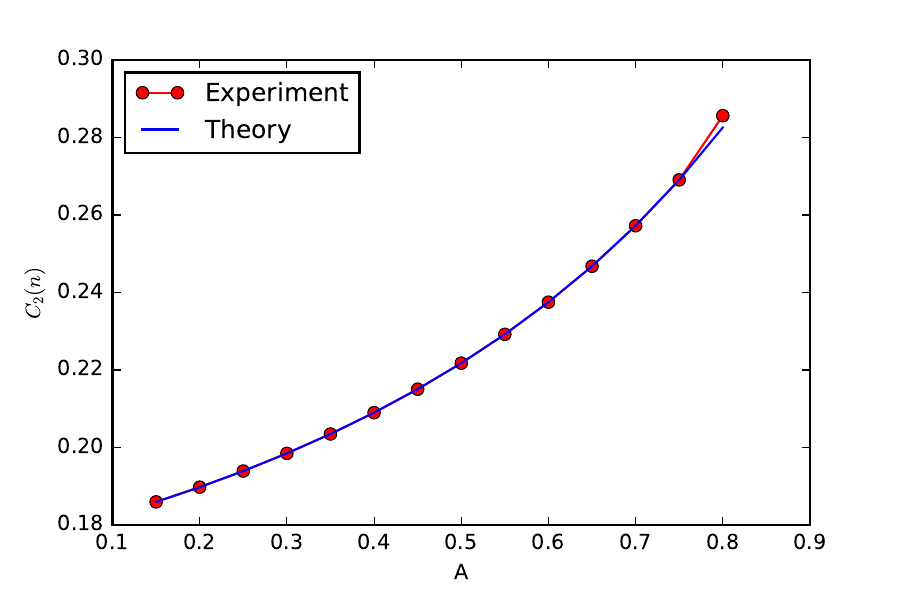}
\end{center}
\caption{The behavior of $C_2(n)$ as a function of $A$ for $n=10^6$, $m=2$, $D=0.3$}
\label{fig:Clust_C_2}
\end{figure}

\textbf{Average Local Clustering Coefficient $C_2(n)$ depending on $A$.}
We generated polynomial graphs with $n = 10^6$, $m = 2$, and $D = 0.3$, assigning $A \in [0.15, 0.8]$. For each value of $A$ we generate 10 graphs and average the obtained values of $C_2(n)$ (see Figure~\ref{fig:Clust_C_2}). For $A \le 0.75$ the theoretical value $\sum_{d=m}^{\infty} {f(d)}$ is extremely close to the experiment and only for $A=0.8$ we observe a small error.
This is consistent with Figure~\ref{fig:C_D_graphs}, where we demonstrated that our approximation of $C(d)$ does not work for $A>\frac 3 4$.

\begin{figure}[t]
\begin{center}
\includegraphics[height = 7cm]{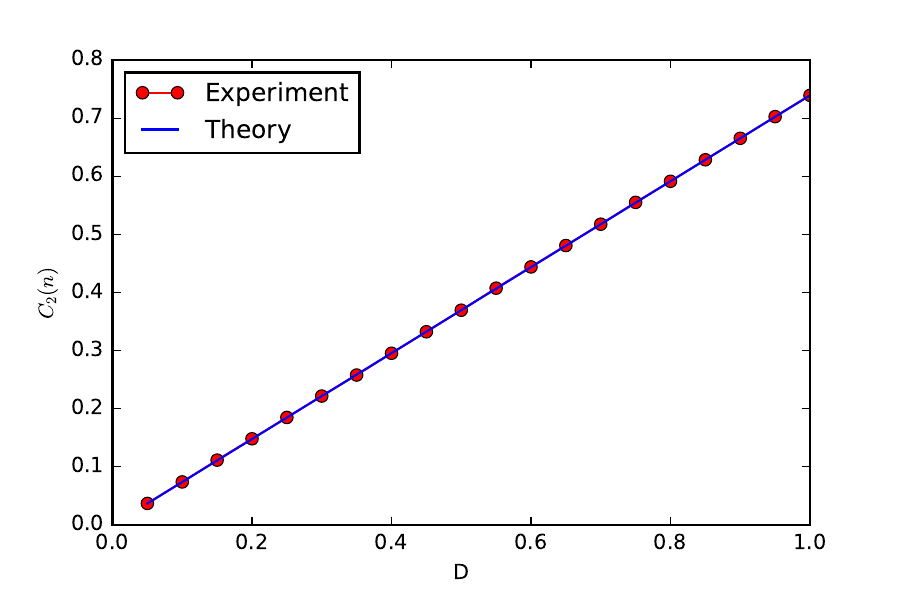}
\end{center}
\caption{The behavior of $C_2(n)$ as a function of $D$ for $n=10^6$, $m=2$, $A=0.5$}
\label{fig:Clust_C_2_D}
\end{figure}

\textbf{Average Local Clustering Coefficient $C_2(n)$ depending on $D$.}
We also generated polynomial graphs with $n = 10^6$, $m = 2$, and $A=0.5$, assigning $D \in [0.05, 1]$. Again, we average $C_2(n)$ over 10 graphs (see Figure~\ref{fig:Clust_C_2_D}).
For all $D$ the theoretical value $\sum_{d=m}^{\infty} {f(d)}$ is extremely close to the experiment.
Also, it follows from Equation~\eqref{C_2D} that $C_2(d)$ should depend linearly on $D$ and our experiment confirmed it.

Thus, our experiments suggest that for polynomial models we can approximate the local clustering coefficient $C_2(n)$ by $\sum_{d=m}^{\infty} {f(d)}$ for $A \le \frac{3}{4}$.

\section{Conclusion}\label{Conclusion}

In this paper, we study the local clustering coefficient $C(d)$ for the vertices of degree $d$ in the T-subclass of the PA-class of models.
Despite the fact that the T-subclass generalizes many different models, we are able to analyze the local clustering coefficient for all these models.
Namely, we proved that $C(d)$ asymptotically decreases as
$\frac{2D}{Am} \cdot d^{-1}$.
In particular, this result implies that one cannot change the exponent $-1$
by varying the parameters $A, D$, and $m$.
This basically means that preferential attachment models in general are not flexible enough to model $C(d) \sim d^{-\psi}$ with $\psi \neq 1$. In addition, we suggested and empirically
verified (for $A \le 0.75$) an approximation for the local clustering coefficient $C_2(n)$.

We would also like to mention the connection between the obtained behavior of $C(d)$ and the notion of \textit{weak} and \textit{strong transitivity} introduced in \cite{Serrano}.
It was shown in  \cite{Serrano2} that percolation properties of a network are defined by the type (weak or strong) of its connectivity.
Interestingly, a model from the T-subclass can belong to either weak or strong transitivity class:
if $2D < Am$, then we obtain the weak transitivity; if  $2D > Am$, then we obtain the strong transitivity.


\newpage

\begin{thebibliography}{30}

\bibitem{BA_Review} R. Albert, A.-L. Barab\'asi, Statistical mechanics of complex networks, Reviews of modern physics, vol. 74, pp. 47--97 (2002)

\bibitem{BioInfoPrior} S. Bansal, S. Khandelwal, L.A. Meyers, Exploring biological network structure with clustered random networks, BMC Bioinformatics, 10:405~(2009)

\bibitem{BA1} A.-L. Barab\'asi, R. Albert, Emergence of scaling in random networks, Science 286, pp. 509-512 (1999)
    
\bibitem{BA2} A.-L. Barab\'asi, R. Albert, H. Jeong, Mean-field theory for scale-free random networks, Physica A 272, pp. 173-187 (1999)


\bibitem{Networks} S. Boccaletti, V. Latora, Y. Moreno, M. Chavez, D.-U. Hwang, Complex networks: Structure and dynamics, Physics reports, vol. 424(45), pp. 175-308 (2006)

\bibitem{Math_Results} B. Bollob\'as, O.M. Riordan, Mathematical results on scale-free random graphs, Handbook of Graphs and Networks: From the Genome to the Internet, pp.1-34~(2003)

\bibitem{LCD_degrees} B. Bollob\'as, O.M. Riordan, J. Spencer, G. Tusn\'ady, The degree sequence of a scale-free random graph process, Random Structures and Algorithms, vol. 18(3), pp. 279-290~(2001)

\bibitem{Chayes} C. Borgs, M. Brautbar, J. Chayes, S. Khanna, B. Lucier, The power of local information in social networks, preprint~(2012)

\bibitem{Broder} A. Broder, R. Kumar, F. Maghoul, P. Raghavan, S. Rajagopalan, R. Stata, A. Tomkins, J. Wiener, Graph structure in the web, Computer Networks, vol. 33(16), pp. 309-320~(2000)

\bibitem{Buckley_Osthus} P.G. Buckley, D. Osthus, Popularity based random graph models leading to a scale-free degree sequence, Discrete Mathematics, vol. 282, pp. 53-63~(2004)

\bibitem{Catanzaro} M. Catanzaro, G. Caldarelli, and L. Pietronero, Phys. Rev. E 70, 037101~(2004)









\bibitem{F-F-F} M. Faloutsos, P. Faloutsos, Ch. Faloutsos, On power-law relationships of the Internet topology, Proc. SIGCOMM'99~(1999)

\bibitem{Girvan} M.~Girvan and M.~E. Newman, Community structure in social and biological networks, Proceedings of the National Academy of Sciences,  99(12):7821-7826~(2002)



\bibitem{Holme_Kim} P. Holme, B.J. Kim, Growing scale-free networks with tunable clustering, Phys. Rev. E, vol. 65(2), 026107~(2002)


\bibitem{Lescovec} J. Leskovec, Dynamics of Large Networks, ProQuest, 2008.

\bibitem{Newman2} M.E.J.~Newman, Power laws, Pareto distributions and Zipf's law, Contemporary Physics, 46, N5, 323-351~(2005)

\bibitem{Newman3} M.E.J.~Newman, The structure and function of complex networks,
SIAM review, 45(2):167-256~(2003)

\bibitem{GPA} L. Ostroumova, A. Ryabchenko, E. Samosvat, Generalized Preferential Attachment: Tunable Power-Law Degree Distribution and Clustering Coefficient, Proc. WAW'13, Lecture Notes in Computer Science, vol. 8305, pp. 185--202 (2013)

\bibitem{Ravasz} E. Ravasz and A.-L. Barab\'{a}si. Hierarchical organization in complex networks. Physical Review E, 67(2)~(2003)



\bibitem{Serrano} M. A. Serrano and M. Bogu\~{n}\'{a},
Clustering in complex networks. I. General formalism,
Phys. Rev. E 74, 056114~(2006)

\bibitem{Serrano2} M. A. Serrano and M. Bogu\~{n}\'{a},
Clustering in complex networks. II. Percolation properties,
Phys. Rev. E 74, 056115~(2006)

\bibitem{Vazquez} A. V\'{a}zquez, R. Pastor-Satorras, and A. Vespignani,
Large-scale topological and dynamical properties of the Internet,
Phys. Rev. E 65, 066130~(2002)


\bibitem{Watts} D.~J. Watts and S.~H. Strogatz, Collective dynamics of 'small-world' networks, Nature 393, pp. 440--442~(1998)

\bibitem{RAN} T. Zhou, G. Yan and B.-H. Wang, Maximal planar networks with large clustering coefficient and power-law degree distribution, Phys. Rev. E, vol. 71(4)~(2005)

\end{thebibliography}
\end{document}